\documentclass[11pt]{amsart}
\usepackage{color}
\usepackage{amsfonts}
\usepackage{amssymb}
\usepackage{amsmath}
\usepackage{amsthm}
\usepackage{amsopn}
\usepackage{graphicx}
\usepackage{times}
\usepackage{url}
\usepackage{MnSymbol}
\usepackage{enumerate}
\usepackage{hyperref}
\hypersetup{
    colorlinks=true,
    citecolor=red,
    filecolor=black,
    linkcolor=black,
    urlcolor=blue
}

\title{New results on torus cube packings and tilings}

\author[M. Dutour Sikiri\'c]{Mathieu Dutour Sikiri\'c}
\address{Rudjer Boskovi\'c Institute, Bijenicka 54, 10000 Zagreb, Croatia}
\email{mdsikir@irb.hr}

\author[Y. Itoh]{Yoshiaki Itoh}
\address{The Institute of Statistical Mathematics, Midori-cho, Tachikawa, Tokyo 190-8562, Japan}
\email{itoh@ism.ac.jp}

\def\QuotS#1#2{\leavevmode\kern-.0em\raise.2ex\hbox{$#1$}\kern-.1em/\kern-.1em\lower.25ex\hbox{$#2$}}

\begin{document}
\newcommand{\RR}{\ensuremath{\mathbb{R}}}
\newcommand{\NN}{\ensuremath{\mathbb{N}}}
\newcommand{\QQ}{\ensuremath{\mathbb{Q}}}
\newcommand{\CC}{\ensuremath{\mathbb{C}}}
\newcommand{\ZZ}{\ensuremath{\mathbb{Z}}}
\newcommand{\TT}{\ensuremath{\mathbb{T}}}
\newcommand{\HH}{\ensuremath{\mathbb{H}}}

\newcommand{\dbcomment}[1]{\marginpar{\footnotesize {\bf Comment:} {#1} }}
\theoremstyle{definition}
\newtheorem{theorem}{Theorem}
\newtheorem{proposition}[theorem]{Proposition}
\newtheorem{corollary}[theorem]{Corollary}
\newtheorem{lemma}[theorem]{Lemma}
\newtheorem{problem}[theorem]{Problem}
\newtheorem{conjecture}[theorem]{Conjecture}
\newtheorem{claim}[theorem]{Claim}
\newtheorem{remark}[theorem]{Remark}
\newtheorem{definition}{Definition}
\newtheorem{example}{Example}
\newtheorem{algorithm}[theorem]{Algorithm}

\begin{abstract}
We consider sequential random packing
of integral translate of cubes $[0,N]^n$ into the
torus $\QuotS{\ZZ^n}{2N\ZZ^n}$.
Two special cases are of special interest:
\begin{enumerate}
\item The case $N=2$ which corresponds to a discrete case of tilings (considered in \cite{cubetiling,book})
\item The case $N=\infty$ corresponds to a case of continuous tilings (considered in \cite{combincubepack,book})
\end{enumerate}
Both cases correspond to some special combinatorial structure and we describe here new developments.
\end{abstract}

\maketitle

\section{introduction}

A {\em cube tiling} is a family $(v_i + [0,1]^n)$ of translates of the unit cube $[0,1]^n$ that tiles $\RR^n$ by translation.
A cube tiling is said to be of {\em class} ${\mathcal T}_N$ if it is $2\ZZ^n$ periodic and if the vectors $v_i$ belong to $\frac{1}{N} \ZZ^n$. Up to scaling this corresponds to $2N\ZZ^n$ periodic tilings of $\RR^n$ by integral translates of the cube $[0,N]^n$.

The origin of the subject of cube packing is with what is called Keller's conjecture (\cite{keller}) that generalize a previous conjecture of Minkowski.
The conjecture states that in every packing of $\RR^n$ by translates of the cube $[0,1]^n$ there exist two cubes that share a facet.
The conjecture was proved to be true for $n\leq 6$ in \cite{perron}. It was proved to be false for $n\geq 10$ in \cite{lagarias} and for $n\geq 8$ in \cite{dim8}. The counterexample found were of class ${\mathcal T}_2$.
Previously, it was found in \cite{szabo} that Keller's conjecture is true for all $n$ if and only if it is true for all cube packings of class ${\mathcal T}_2$ and all dimensions.

If one restricts to the case of cube tilings of class ${\mathcal T}_2$ of dimension $n$ then the Keller's conjecture is a finite problem (\cite{corradi}).
It is equivalent to proving that some graph $G_n$ having $4^n$ vertices has a clique number lower than $2^n$. In \cite{Dim7case} it was proven that the clique number of this graph $G_7$ is $124$. This shows that Keller's conjecture is true in dimension $7$ for the class of ${\mathcal T}_2$ cube tilings. But it does not a priori rule out the possibility of a counter-example in dimension $7$ that is not of class ${\mathcal T}_2$, though this is unlikely. The computation in \cite{Dim7case} was an extraordinary accomplishment that was for long though to require a computer, the ``size of a galaxy''.

It therefore appears that the class ${\mathcal T}_2$ of cube tilings is very interesting to study. Also a priori interesting is to consider cube packings. For such classes ${\mathcal P}_N$, questions of iterative packing and extensibility occurs.
Also, since the class ${\mathcal T}_2$ is finite in a given dimension, questions of classifications occurs that may be of interest. We will report on the recent works in the subject.

The classes ${\mathcal T}_N$ for $N>2$ are much harder to study combinatorially and of limited interest. However, as $N$ goes to $\infty$, one can study a special kind of continuous cube tiling obtained with positive probability.
This was introduced in \cite{combincubepack} but there has been less progress on this problem since it is less directly combinatorial. We will therefore instead report on the problems that we consider the most important.

\section{Cube packings and tilings for $N=2$}

A {\em cube packing} is a family $(v_i + [0,1]^n)_{i\in I}$ of translates of the unit cube $[0,1]^n$ that tiles $\RR^n$ by translation.
A cube packing is said to be of {\em class} ${\mathcal P}_N$ if it is $2\ZZ^n$ periodic and if the vectors $v_i$ belong to $\frac{1}{N} \ZZ^n$.

One method for obtaining cube tilings of class ${\mathcal T}_2$ is to take vectors $v$ at random in $\left\{0, \frac{1}{N}, \dots, \frac{2N-1}{N}\right\}^n$ and add the cubes if they do not overlap with preexisting ones. This method can be considered as a random process and this approach has a long history starting from \cite{renyi1958}. It was extended to packing in the cube $[0,4]^n$ by integral translates of the cube $[0,2]^n$ in \cite{ueda,dip,poyarkov2} where estimates on the expectation of the obtained random cube packing are obtained \cite{book}.

\begin{figure}
\begin{center}
\begin{minipage}{5.2cm}
\centering
\resizebox{5.0cm}{!}{\includegraphics[bb=105 260 494 559, clip]{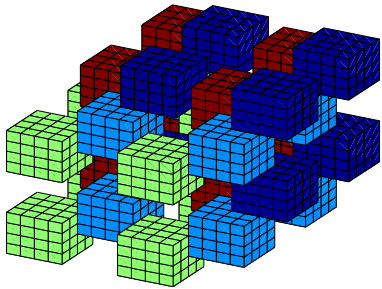}}\par
\end{minipage}
\end{center}
\caption{The unique non-tiling non-extensible cube packing in dimension $3$}
\label{UniqueInex}
\end{figure}

In dimension $n\leq 2$ the sequential random cube packing into torus will always give a tiling. However, in dimension $3$ we can obtain with non-zero probability a non-extensible cube packing with $4$ translation classes and thus density $1/2$ (see Figure \ref{UniqueInex}).

The full classification of ${\mathcal T}_2$ cube tilings has been obtained in dimension $n\leq 5$. See Table \ref{TableClassification} for the number of types.
Since for $n\leq 5$, Keller's conjecture is true any cube tiling has two cubes that share a face. This means that we can shift them by $\frac{1}{2}e_i$ with $e_i$ a basis vector. By iterating such operations one can get new cube tilings. It turns out that for $n\leq 5$ all cube tilings can be obtained by starting from the trivial cube tiling by integral translates of $[0,1]^n$ (see \cite{CubeTilingDim5} for the case $n=5$ and \cite{cubetiling} for $n\leq 4$). It would be interesting to know if this also holds in dimension $6$ and $7$.

The packings of low density are harder to study. Let us denote by $f(n)$ the minimum number of translation in a non-extensible cube packing of class ${\mathcal P}_2$. Similarly let us denote by $h(n)$ the minimum number of translation classes of cubes, possibly overlapping that are needed in order to prevent the addition of one non-overlapping cube (see \cite{cubetiling} for details). Obviously $f(n)\geq h(n)$.
Table \ref{KnownValues_F_H} from \cite{BlockingSet} gives the known values of $f(n)$ and $h(n)$. Other exhaustive enumeration works seems impossible in that direction now. However, one direction that has not been really considered is obtaining infinite families of low density packings for all $n$.

A {\em hole} $H(P)$ is the complement of a non-extensible cube packing $P$ with $2^n - l$ cubes.
Such non-extensible cube-packings do not exist for $l=1$, $2$ or $3$ and any dimension $n$ (\cite[Theorem 2]{cubetiling}).
In \cite[Conjecture 1]{cubetiling} we stated that for $l=4$ the holes $H(P)$ is essentially unique and is given by the one of Figure \ref{UniqueInex} and its higher dimensional extensions.
We also conjectured that holes do not exist for $l=5$. The conjecture was proved in \cite{EnumerationHoles} for $n=5$ by an exhaustive enumeration.
Also we conjectured that for $l=6$ or $l=7$ in any dimension, the holes belong to a finite set of possibilities.
It seems to us that such conjectures are not hopeless and could be proved by extending the proof technique of \cite[Theorem 2]{cubetiling}.

\begin{table}
\caption{Number of types of cube tilings in dimension $n\leq 5$}
\label{TableClassification} 
\begin{center}
\begin{tabular}{|c|c|c|c|c|}
\hline
n       & $\#$ types        & Reference\\
\hline
$2$     &   $2$             &\\
$3$     &   $9$             &\\
$4$     &   $744$           & \cite{cubetiling}\\
$5$     &   $899,710,227$   & \cite{CubeTilingDim5}\\\hline
\end{tabular}
\end{center}
\end{table}

\begin{table}
\caption{Known values and ranges of $f(n)$ and $h(n)$ for $n\leq 7$ (from \cite{BlockingSet})}
\label{KnownValues_F_H}
\begin{center}
\begin{tabular}{|c|c|c|}
\hline
n       &  $f(n)$           & $h(n)$\\
\hline
$2$     &   $4$             &$3$\\
$3$     &   $4$             &$4$\\
$4$     &   $8$\cite{cubetiling}      &$7$\cite{cubetiling}\\
$5$     &   $12$\cite{BlockingSet}      &$10$\cite{brink}\\
$6$     &   $16$\cite{BlockingSet}      &$15$\cite{BlockingSet}\\
$7$     &   $20-32$\cite{BlockingSet}      &$20-23$\cite{BlockingSet}\\\hline
\end{tabular}
\end{center}
\end{table}

\section{Continuous cube tilings and packings}

The main peculiarity of cube tilings of class ${\mathcal T}_N$ is the $2\ZZ^n$ periodicity. It imposes that for any two translation classes of cubes $t + [0,1]^n$ and $t' + [0,1]^n$ that are non-overlapping, there exist a coordinate $i$ such that $t_i$ and $t'_i$ differ by an integer.
We choose to consider this as their main feature and see the behavior of the cube tilings as $N\to \infty$.
The formalism is explained in details in \cite{combincubepack} and allows to consider the cube packing that are obtained with positive probability as the main objects.

From the viewpoint of exhaustive combinatorial enumeration the problems are easier with a slower combinatorial explosion (for example in dimension $4$ we have $32$ types of such continuous cube tilings vs $744$ for the class ${\mathcal T}_2$). However, the absence of a graph formalism makes it harder to program and they were thus much less studied.

On the other hand the continuous structure gives the notion of number of parameters that are needed to describe the structure. It is conjectured (\cite[Conjecture 5.4]{combincubepack}) that this number of parameters is at most $2^n -1$ but we were unable to prove it. The number of parameters is at least $\frac{n(n+1)}{2}$ but an open question is to prove the existence of a cube tiling with this number of parameters and obtained with positive probability.

In one respect the continuous case is simpler. For $n$ odd the minimal non-extensible cube packings can be classified (\cite[Proposition 5.5]{combincubepack}): they have $n+1$ translation classes of cubes and are described by what are called one-factorizations of perfect graphs $K_{2m}$ with $2m=n+1$. Such one-factorizations exist for any $m\geq 1$ (\cite{walecki}) and the number of non-isomorphic types is known for $m\leq 7$ (\cite{KO_K14}).
If $n$ is even, then one expects the existence of a non-extensible continuous cube packing with $n+2$ cubes and $\frac{n(n+1)}{2}$ parameters, but the expected cube packings would be more complicated than one-factorization.

\section*{Acknowledgments}
Both authors thank the Institute of Statistical Mathematics for support.

\end{document}